\documentclass[letterpaper, 10pt, conference]{ieeeconf}
\IEEEoverridecommandlockouts

\usepackage{graphicx}
\usepackage{float}
\usepackage{subfigure}
\usepackage{derivative}
\usepackage[utf8]{inputenc}
\usepackage{nicefrac}

\usepackage{bm}
\usepackage{mathptmx}
\usepackage{tcolorbox}
\usepackage{tikz}
\usepackage{amssymb, amsmath}
\usepackage{mdframed}

\usepackage{dsfont}
\usepackage{graphicx}
\usepackage[colorlinks=true, allcolors=blue]{hyperref}
\newcommand{\p}{{\mathbb P}}

\newcommand{\pb}{{Q}}
\newcommand{\bb}{{\bm b}}
\newcommand{\A}{{\cal A}}
\newcommand{\s}{{\cal S}}

\newcommand{\E}{{\mathbb{E}}}

\newcommand{\eop}{{\hfill $\Box$}}
\newtheorem{thm}{{\bf Theorem}}
\newtheorem{lem}{{\bf Lemma}}

\newcommand{\ignore}[1]{}

\IEEEoverridecommandlockouts                              

\usepackage{bbm}
\usepackage{enumerate}

\newcommand{\exactproof}[1]{}  

\begin{document}
\title{\bf Constrained Average-Reward   Intermittently Observable MDPs \\
}

\author{Konstantin Avrachenkov, Madhu Dhiman and Veeraruna Kavitha 
\thanks{K.~Avrachenkov is with INRIA Sophia Antipolis, France. Email: k.avrachenkov@inria.fr; M.~Dhiman and V.~Kavitha are with the Indian Institute of Technology Bombay, India. Emails: \{madhu.dhiman,vkavitha\}@iitb.ac.in. The authors are listed in alphabetical order. This work is partially supported by DST-INRIA Cefipra project  Learning In Operations and Networks (LION). The paper is accepted at IEEE CDC 2025.}
}
\maketitle
\begin{abstract}
In Markov Decision Processes (MDPs) with intermittent state information, decision-making becomes challenging due to periods of missing observations. Linear programming (LP) methods can play a crucial role in solving MDPs, in particular, with constraints. However, the resultant belief MDPs lead to infinite dimensional LPs, even when the original MDP is with  finite state and action spaces. The verification  of strong duality  becomes non-trivial.  
This paper investigates the conditions  for  no duality gap in  average-reward  finite Markov decision process with intermittent state observations. We first establish that in such MDPs, the belief MDP is unichain  if the original Markov chain is recurrent. Furthermore, we establish  strong duality of the problem, under the same assumption. Finally, we provide a wireless channel example, where the belief state depends on the last channel state received and the age of the channel state.  Our numerical results indicate interesting properties of the solution.   
\end{abstract}
\section{Introduction}
\label{sec_intro}

Markov Decision Processes (MDPs) provide a rigorous mathematical framework for sequential decision-making under uncertainty. They have been widely applied across various domains, including wireless communication, finance, and supply chain management. Classical MDPs, as described in \cite{puterman2014markov}, operate under the fundamental assumption that the decision-maker has full knowledge of the system’s state at every decision epoch. However, in many real-world scenarios, this assumption does not hold, as decision-makers often have access to only partial or delayed state information. This challenge gives rise to Partially Observable MDPs (POMDPs), which have been extensively studied, see e.g. \cite{krishnamurthy2016partially} and many references therein. Most of the studies have concentrated on the conditions of the existence of the optimal policies  and dynamic programming principle, \cite{borkar2003dynamic,borkar2007dynamic,demirci2024average,hsu2006existence,kara2022near,platzman1980optimal,runggaldier1993nearly}. However, the   literature that studies the    infinite-dimensional LPs associated with POMDPs  is extremely sparse;  the LP approach for POMDP with constraints has received even lesser attention. 

In many practical applications, such as those discussed in \cite{krishnamurthy2016partially, chen2023intermittently, khan2023cooperative,shuman2010measurement}, state information is not only partial but also intermittently available due to factors like communication constraints, sensor failures, or energy limitations. This motivates the study of Intermittently Observable MDPs (IOMDPs) \cite{chen2023intermittently}, where decision-making must proceed despite missing state observations at certain time steps.  

A common approach to handling POMDPs is to reformulate them as belief MDPs, where the decision-maker maintains a probability distribution (or belief) over the system  states based on past observations and prior knowledge \cite{krishnamurthy2016partially,bertsekas1996stochastic}. While this transformation enables decision-making under uncertainty, it significantly increases computational complexity. Unlike classical MDPs with a finite state space, the belief MDP generally has an uncountable state space, calling for various approximation and state aggregation approaches. Furthermore, even in the case of finite state space models, the MDPs with constraints cannot be solved by value iteration or policy iteration. This challenge has led to the development of Linear Programming (LP)-based methods, which have proven effective for solving constrained MDPs \cite{altman2021constrained}. However, extending LP formulations to IOMDPs is non-trivial, as it requires verifying strong duality, a property that is difficult to establish due to the infinite-dimensional nature of the belief space \cite{lasserre1994average, altman2021constrained}.
\ignore{\color{red}making standard solution techniques such as value iteration and policy iteration impractical. This challenge has led to the exploration of Linear Programming (LP)-based methods, which have proven effective for solving constrained MDPs \cite{altman2021constrained}. However, extending LP formulations to IOMDPs is non-trivial, as it requires verifying strong duality, a property that is difficult to establish due to the infinite-dimensional nature of the belief space \cite{lasserre1994average, altman2021constrained}.}%
\ignore{\color{blue}The  LP based methods   (well-known for solving the problems  with constraints) have similar kind of complexity, but  more importantly face a more fundamental complication. It is non-trivial to  verify strong duality condition  for IOMDPs  due to the infinite-dimensional nature of the belief space \cite{lasserre1994average, altman2021constrained}.} In particular, the results in \cite{lasserre1994average} require  a certain set of measures being tight (see Assumption~A.3 there). However, when the cost function is bounded, as considered in this work, the relevant set of measures becomes significantly large, which is not tight with respect to the  required weak topology in \cite{lasserre1994average}; and thus one cannot directly apply the results of \cite{lasserre1994average} to our framework.

Ensuring strong duality is crucial, as it facilitates the efficient computation of optimal policies and extends LP methods to constrained IOMDPs, which are relevant in applications requiring adherence to performance guarantees or safety constraints.  

Most of the existing literature on POMDPs and IOMDPs \cite{chen2023intermittently, khan2023cooperative, lim2023optimality, ahmadi2020control} focuses on the discounted cost/reward criterion. The discounted objective has several analytical advantages, including contraction properties that simplify the analysis and guarantee the existence of optimal policies. In \cite{chen2023intermittently}, the authors analyze an IOMDP under a total expected discounted reward criterion. Our work is inspired by this study but differs in that we consider the average-reward criterion, which presents additional challenges due to the absence of discounting-related contraction properties. The POMDP formulation in general leads to a belief MDP with an uncountable state space. However, following the insights from \cite{chen2023intermittently}, we demonstrate that  a countable belief state representation is sufficient in IOMDPs. In particular, this allows us to handle belief MDPs with a convenient approximation, making the problem more tractable.  

Specifically, in this paper, we investigate the conditions under which strong duality holds in an average-reward IOMDP with finite state space and intermittent state observations. Assuming that the underlying Markov chain is recurrent, we first establish that the corresponding belief MDP is unichain. We then prove that the strong duality holds under the same assumption.  This establishes that both the constrained and unconstrained LP formulations of `recurrent' IOMDPs do not exhibit  duality gap. Thus, our results provide a solid foundation for applying LP techniques to solve IOMDPs, both in constrained and unconstrained cases. 

Finally, we apply our theoretical results to obtain an (approximate) optimal policy for an application with constraints,  that observes the state  of the system over a lossy wireless channel. By solving the approximate LPs (obtained after excluding the `large' age states), we  compute the optimal policy that depends upon the last observed state and the age of the information; we further discuss the properties of the solution.

\vspace{-2mm}
\section{ Constrained Finite MDP Model}
\label{sec_pre}
A discrete time constrained Markov decision process (MDP) is represented by the tuple $({\cal S}, {\cal A}, {P}, r, c, B)$ where ${\cal S}$ is a  set of states, ${\cal A}$ is a set of actions, $r$ is an immediate reward function, $c$ is an immediate cost function, $B$ is a finite bound on the long run average cost (see \eqref{Eqn_average_reward_cons_underline}) and  $P$ is a controlled transition kernel $P : {\cal S} \times {\cal A} \times {\cal S} \rightarrow[0,1]$.

The transition kernel
$P(s'|s,a)$ represents the probability of arriving in state $s^{\prime}$ having taken action $a$ in state $s$:
$$
P_a(s, s')=P(s' \mid s, a)=\operatorname{Pr}\left(s_{t+1}=s' \mid s_t=s, a_t=a\right), \ \forall t,
$$ 
where  $s_t$ and $a_t$ denote the state and action at time $t$.  
%
We will use $P_a=\left[P_a(s, s')\right]_{s, s'}$ to denote the transition matrix associated with action $a$. Note that without loss of generality we are assuming the same action space for each state.

The immediate reward function $r: {\cal S} \times {\cal A} \rightarrow \mathbb{R}$ (set of real numbers), with  $r(s,a)$ representing the immediate expected reward of taking action $a$ in state $s$. In a number of applications, it is useful to consider also the immediate cost function $c: {\cal S} \times {\cal A} \rightarrow \mathbb{R}$.  In this work, we assume finiteness of the state space and the action space:  

\medskip
\noindent\textbf{A.1} \textit{The state and action spaces, 
${\cal S}$ and ${\cal A}$ are   finite}.
\medskip

We also assume the following recurrence condition, which is satisfied in many practical cases:  

\medskip
\noindent \textbf{A.2} \textit{The  MDP $({\cal S}, {\cal A}, {P}, r, c, B)$ is recurrent}.
\medskip

Consider a set of stationary Markov and randomized  (SMR) policies  for such MDPs \cite{puterman2014markov,altman2021constrained}; any SMR policy is denoted by $\pi: {\cal S} \rightarrow \p ({\cal A})$,   a function that maps the set of states to the set of probability measures on ${\cal A}$. 

The problem of a constrained average-reward MDP can be formulated as follows:
\begin{equation}
    \sup_{\pi \in  \mbox{SMR}}\left(\lim_{T \to \infty} \frac{1}{T} \sum_{t=1}^T \E[ r(s_t, a_t)]  \right),
\label{Eqn_average_reward_uncons_underline}
\end{equation}
subject to
\begin{equation}
\label{Eqn_average_reward_cons_underline}
 \lim_{T \to \infty} \frac{1}{T} \sum_{t=1}^T \E[ c( s_t, a_t )]  \le B.
\end{equation}
Under assumptions \textbf{A.1-2}, the above problem has a solution in the set of SMR policies (see e.g., \cite{altman2021constrained}).


\section{IOMDP Model}
\label{sec_model}
In classical MDP, the system executes the optimal policy $\pi$ in the following manner: at every time step, it observes its current state $s_t$, and chooses an action according to the probability measure $\pi(s_t)$. However, in some applications, the system may not be able to observe current state with perfect information (e.g., states observed over wireless medium are prone to erasures and/or errors, estimation methods are prone to errors, less expensive to observe limited information, etc.).  In such 
partially observable MDPs, the system may observe the current state with varying degrees of information; it can observe some partial information about the current state or may  observe none.  In this paper, we consider a  special case where the system   either observes the  state completely with a fixed probability  or receives zero information. Such systems are called Intermittently Observable MDPs (IOMDPs).  
In  \cite{chen2023intermittently}, the authors consider  an IOMDP with discounted rewards, and we consider a similar setup but for average-reward criteria; we further  consider  a constrained problem. 

Namely, at every time step, the controller observes the current state completely with probability $\rho \in (0, 1]$ or does not observe with complementary probability $1-\rho$.  
The observation events are independent and identically distributed across the time steps.


\subsection{Belief MDP}
As is usually done in POMDP literature (see e.g., \cite{krishnamurthy2016partially}), we reformulate the problem as a belief MDP. Towards this, we  consider `belief states',  the conditional distribution of the current state, given the history of observations.


 The belief state at time 
$t$, denoted by $\bb_t$, represents  the probabilistic estimate of the system state at the same time, given all the available information until $t$; this forms a  sufficient statistic for the problem (see \cite{krishnamurthy2016partially,bertsekas1996stochastic}).
Let $\ell_t$ represent the last observation  epoch before time $t$. Then, the belief state $\mathbf{b}_t$, equals:
\begin{eqnarray}
    \bb_t & :=& \operatorname{Pr}\left(s_t \mid s_{\ell_t}, a_{\ell_t}, a_{\ell_t+1}, \cdots, a_{t-1}\right).
    \label{eqn_belief_state}
\end{eqnarray}
Basically, by Markovian nature of the problem and under SMR policies, the belief state is just  the function of the last observed state $s_{\ell_t}$ and the sequence of actions thereafter, $(a_{\ell_t}, a_{\ell_t+1}, \cdots, a_{t-1})$. 
In other words, the sequence $\left(s_{\ell_t}, a_{\ell_t}, a_{\ell_t+1}, \cdots, a_{t-1}\right)$ represents the sufficient statistic at time $t$. 

 If state $s_t$ is observed at time $t$ (and say $s_t = i$), then the belief state $\bb_t$  is updated to the standard basis vector $\mathbf{e}_i$ (vector with $i$-th component as one and all other components being zero). And, if state $s_t$ is not observed (which happens with probability $(1-\rho)$), then the  belief state $\bb_t$ is updated based on the  belief $\bb_{t-1}$ and the action $a_{t-1}$ of the previous epoch.   Let $\pb$ represent the transition kernel of the belief states. Using the above arguments,  $Q$ is given~by
\begin{equation}
\label{eqn_belief_trans_probability}
\pb \left(\bb_t \mid \bb_{t-1}, a_{t-1}\right) :=   \left\{\begin{array}{ll}
1-\rho, & \text {if } \bb_t=P_{a_{t-1}}^{\top} \bb_{t-1}  \\
\rho\left[P_{a_{t-1}}^{\top} \bb_{t-1}\right]_i, & \text {if } \bb_t=\mathbf{e}_i, i \in \mathcal{S}  \\
0, & \text {otherwise}.
\end{array}\right.  
\end{equation}
In the above,  $P_{a_{t-1}}^{\top} \bb_{t-1}$
 is  the updated belief distribution, and  $\left[P_{a_{t-1}}^{\top} \bb_{t-1}\right]_i$ denotes the $i$-th element of vector $P_{a_{t-1}}^{\top} \bb_{t-1}$; here~$(\cdot)^\top$ represents the transpose. We note that the evolution of the belief state \eqref{eqn_belief_trans_probability} is similar to the dynamics considered in \cite{shuman2010measurement,akbarzadeh2024two}. Let $\mathcal{B}$ represent the set of belief states. From \eqref{eqn_belief_trans_probability},  we have: 
\begin{equation*}
    \mathcal{B} := \left\{\prod_{k=0}^L P_{a_k}^{\top} \mathbf{b}: a_k \in \mathcal{A}, \forall k \le L, \ \mathbf{b} \in {\mathcal{E}} \mbox{ and }L \in \{ 0,1,  \cdots\}\right\},
\end{equation*}
where  $\mathcal{E}:= \left\{ \mathbf{e}_i: i \in \mathcal{S}\right\}$ represents the set of pure (or just observed) states.
%
 Define the function $R: \mathcal{B} \times \mathcal{A} \to \mathbb{R}$ as below:
\begin{equation}
    R(\bb, a) = \sum_{s \in \mathcal{S}} \bb(s) r(s,a), \ \forall \bb \in \mathcal{B}, \ a \in \mathcal{A}. 
    \label{Eqn_belief_reward}
\end{equation}
Observe here that, for IOMDPs, under any SMR  policy, the  actions chosen  depend (only) upon the belief states at the respective decision epochs.
Let $\Pi^{b}$ represent the set of such SMR policies. Consider any $\pi \in \Pi^b$ and then using~\eqref{Eqn_belief_reward},
\begin{equation}
     \lim_{T \to \infty} \frac{1}{T} \sum_{t=1}^T \E[ R(\bb_t, a_t)]  = \lim_{T \to \infty} \frac{1}{T} \sum_{t=1}^T  \E \left [  \sum_{s \in \mathcal{S}} \bb_t(s) r(s,a_t) \right ].  \label{Eqn_average_reward_uncons_belief}
\end{equation} 
By A.1-2, the  limit in the above exists\footnote{Any SMR policy based on belief states is a history dependent policy of the original MDP and by \cite[Theorem 6.9.1]{puterman2014markov}, see also \cite{altman2021constrained}, the performance under each such policy is achieved by that under an SMR policy and hence, the limit exists by A.1-2. Further, the problem in \eqref{Eqn_average_reward_uncons_belief} can be seen as a problem under certain `observability' constraints (see~\eqref{Eqn_belief_reward}).}, as it  equals the   average reward of the original MDP under an appropriate history dependent policy  $\pi$.  
Thus,  $R$ can be viewed as  the immediate reward function defined in terms of  the belief states. 
Further, boundedness of  $r(\cdot, \cdot)$ implies that of   $R(\cdot, \cdot)$. 
Finally,  the belief MDP is given by tuple $(\mathcal{B}, \mathcal{A}, \pb, \rho, R)$. 
Now, the unconstrained problem for the IOMDP can be defined by the following: 
\begin{equation}
    \sup_{\pi \in \Pi^b}\left(\lim_{T \to \infty} \frac{1}{T} \sum_{t=1}^T \E[ R(\bb_t, a_t)]  \right).
\label{Eqn_average_reward_uncons_belief}
\end{equation}
By Theorem \ref{Thm_two_parts} given below, the belief MDP is unichain under A.1-2 and hence the above problem  is well defined.
Likewise, we also consider  the constrained (see \eqref{Eqn_average_reward_cons_underline}) average reward belief  MDP, where  we aim to solve:
\begin{eqnarray}
\label{Eqn_average_reward_cons_belief}
\sup_{\pi \in \Pi^b} \hspace{-3mm}&\hspace{-3mm}&\hspace{-3mm} \left(\lim_{T \to \infty} \frac{1}{T} \sum_{t=1}^T \E[ R(\bb_t, a_t)]  \right), 
\mbox{ subject to }  \nonumber \\
 && \hspace{-3mm} \lim_{T \to \infty} \frac{1}{T} \sum_{t=1}^T \E[ C( \bb_t, a_t )]  \le B, \mbox{ where}  \\
  C(\bb, a) &:=& \sum_{s \in \mathcal{S}} \bb(s) c(s,a), \ \forall \bb \in \mathcal{B}, \ a \in \mathcal{A}. \nonumber 
\end{eqnarray}
In the above,   $C (\cdot,\cdot)$ is defined similarly to \eqref{Eqn_belief_reward}. 
We now proceed to proving that the belief MDP is unichain. Towards that, we require some definitions of \cite{ altman2021constrained}, which we reproduce here in our notations:
\begin{itemize}
\item \textbf{$\mu$-bounded function space, $F^{\mu}$:} \textit{For given $\mu: \mathcal{B} \to \mathbb{R} $, $F^{\mu}$ is space of all $\mu$-bounded functions, i.e., 
$$
F^{\mu}:= \left\{q: \mathcal{B} \to \mathbb{R}; \ \sup_{\bb \in \mathcal{B}} \left(\frac{q(\bb)}{\mu(\bb)} \right)< \infty \right\}.
$$}
    \item  \textbf{Contracting MDP:} \textit{Let ${\mathcal{B}}^{\prime}$ and $\mathcal{E}$ be two disjoint sets of states with ${\mathcal{B}}={\mathcal{B}}^{\prime} \cup \mathcal{E}$. An MDP is said to be contracting (on ${\mathcal{B}}^{\prime}$) if there exist a scalar $\xi \in[0,1$) (called the contracting factor) and a function $\mu: \mathcal{B} \rightarrow[1, \infty)$, such that for all $\bb \in \mathcal{B}, a \in \mathcal{A}$,
$$
\sum_{\bb' \notin \mathcal{E}} \pb(\bb'|\bb,a) \mu(\bb') \leq \xi \mu(\bb), 
$$ 
with additional three conditions on  MDP setup:
\begin{enumerate}
    \item $\langle \beta, \mu \rangle = \int \mu(\bb) \beta(d\bb) < \infty$, where $\beta$ is the initial state distribution;
    \item The transition probabilities $Q$ are $\mu$-continuous, i.e., if $a_n \to a$, then  for all $\bb \in \mathcal{B}$
    $$
\lim_{n \to \infty} \sum_{\bb' \in \mathcal{B}} |\pb (\bb'| \bb, a_n)- \pb (\bb'| \bb, a)| \mu(\bb') =0;
    $$
    \item $R(\cdot, \cdot)$ is $\mu$-bounded by some finite constant. Additionally,  for constrained IOMDP problem, $C(\cdot, \cdot)$ is also $\mu$-bounded by some finite constant. 
\end{enumerate}
}
\item \textbf{Uniform Lyapunov function:} \textit{A function $\mu: \mathcal{B} \to [1, \infty) $ is called a uniform Lyapunov function  if there exist constants $ \lambda \in (0,1) $, $ B > 0 $ and a finite set $\mathcal{C} \subset \mathcal{B}$ such that  for all states $ \bb \in \mathcal{B} \setminus \mathcal{C}$ and all $ a \in \mathcal{A}$, the expected drift condition holds:  
$$
\sum_{\bb' \in \mathcal{B} } \pb(\bb'|\bb,a) \mu(\bb') \leq \mu(\bb)  - \lambda.
$$  
}
\end{itemize}

It is clear from \cite[Section 8.2]{puterman2014markov}), every SMR policy $\pi$ will generate a belief Markov chain. Now, we show that the belief MDP is unichain and characterize a few properties, which are useful to prove the ensuing results.
\begin{thm} {\bf [Unichain]}
\label{Thm_two_parts}
    Assume A.1-2, then the following are true: 
    \begin{enumerate}[(i)]
        \item For any SMR policy $\pi$, the belief Markov chain has a single positive recurrent class and absorption into the positive recurrent class  takes place in a finite expected time. 
        \item Belief MDP is contracting  and has a uniform Lyapunov function.
    \end{enumerate}
\end{thm}
\textbf{Proof} is in Appendix \ref{sec_app_a}. \eop

Indeed, by part $(i)$, the belief MDP is  unichain, which means that \eqref{Eqn_average_reward_uncons_belief} and \eqref{Eqn_average_reward_cons_belief}  are well defined. And, part $(ii)$ implies that the belief MDP has a uniform Lyapunov function, which will help in proving the strong duality in the following section.

\vspace{-4mm}
\section{Strong Duality}
\label{sec_strong_duality}
It is well known that one can solve  perfect information finite   MDPs, using linear programming (LP) techniques. 
In fact, the  LP based techniques are extremely useful for solving constrained  finite MDPs.  The quest in this section is to investigate if the same is possible for the case with intermittent observations.  

In the case of IOMDPs, the belief state space is countable. Hence, a~priori, it is not evident if the LP based approach can be used here; the primary concern is related to the duality gap.   The main aim of the paper is to address this issue. We begin with defining the potential LP for unconstrained IOMDPs  in the immediate next.

\subsection{Unconstrained IOMDPs}
\label{sec_uncons_iomdp}
The primal  linear program (LP) for the belief MDP with objective function as in \eqref{Eqn_average_reward_uncons_belief} is given by:
\begin{eqnarray}
\label{Eqn_primal_lp_uncons}    
 \inf_{ {\bf x} = \{x(\bb, a)\}_{\bb, a} }  && \sum_{\bb,a} - R(\bb,a)x(\bb,a) \ \  \mbox{subject to } \nonumber \\
&& \sum_a x(\bb',a) = \sum_{\bb,a} \pb(\bb' | \bb,a) x(\bb,a) , \ \forall \bb' \in \mathcal{B}  \nonumber \\
&&\sum_{\bb,a} x(\bb,a) = 1, \  x(\bb,a) \ge 0.
\end{eqnarray}
 Since our work mostly depends upon the results from \cite{altman2021constrained}, we consider now onwards the cost (which is negative of reward) minimization as in \cite{altman2021constrained}.   
Here, each $\{x(\bb,a)\}_{\bb,a}$ vector is supposed to represent the  occupancy measures for some policy $\pi \in \Pi^b$, where  
\begin{equation*}
    x(\bb,a) := \lim_{T \to \infty} \frac{1}{T}  \sum_{t=1}^T \operatorname{Pr} (\bb_{t} =\bb, a_t =a), \  \forall \bb \in \mathcal{B}, a \in \A. 
\end{equation*}
Since $\bb' \in \mathcal{B}$ (a countable space), the linear program \eqref{Eqn_primal_lp_uncons} operates over an infinite-dimensional space, making the solution procedure non-trivial.

Before stating the result, we  define the standard average cost optimality equation (ACOE). For that, first we define the gain and  the bias functions for the belief MDP. The gain of the Markov chain under policy $\pi$ with the initial state $\bb$ is defined by
$$
g^{\pi}(\bb) :=  \lim_{T \to \infty} \frac{1}{T}\E_{\bb}\left[\sum_{t=1}^{T}   -R\left(\bb_t, \pi(\bb_t)\right) \right], \quad \forall \bb \in \mathcal{B}.
$$
Since $R(\cdot, \cdot)$ is a bounded function and $\mathcal{S}, \A$ are finite,  $g^{\pi}$ is well defined and is also bounded, for every SMR policy $\pi$; also observe $g^{\pi}(\bb)  = g^{\pi}$,  a constant for all belief states $\bb$,   by unichain property proved  in Theorem \ref{Thm_two_parts}. Thus  the bias function $h^{\pi}(\cdot)$ is also well defined and is given by:
\begin{eqnarray}
h^\pi (\bb) &:=& \lim_{T \to \infty} \frac{1}{T } \sum_{k=1}^{T}\mathbb{E}_\bb \left[ \sum_{t=1}^{k}\bigg(- R(\bb_t, \pi(\bb_t))-g^{\pi}(\bb_t)\bigg)\right],  \nonumber  \\ && \hspace{48mm} \forall \bb \in \mathcal{B}. \ \ 
\label{eqn_bias_uncons}
\end{eqnarray}
Using the unichain property proved in Theorem~\ref{Thm_two_parts}.(i),  the ACOE can be written  as (see e.g., \cite{puterman2014markov}) 
\begin{equation}
    \phi(\bb) + \psi = \min_{a \in \mathcal{A}} \left( - R(\bb, a) + \sum_{\bb' \in \mathcal{B}} \pb (\bb'| \bb, a) \phi(\bb')\right), 
    \label{eqn_ACOE}
\end{equation}%
where $\psi$ is some scalar and $\phi: \mathcal{B} \to \mathbb{R}$. If for some SMR policy $\pi$,  equation \eqref{eqn_ACOE} is satisfied 
with $\psi = g^\pi$ and $\phi(\bb) = h^\pi(\bb)$ for all $\bb$,  then  $g^{\pi}$ is the optimal value (see  e.g., \cite{puterman2014markov, altman2021constrained}). Now, the  dual problem  of \eqref{Eqn_primal_lp_uncons} is given by   
\begin{eqnarray}
\label{Eqn_dual_lp_uncons}    
 \sup_{\psi, \phi} && \psi \quad  \mbox{subject to }  \\
&& \phi(\bb)+ \psi \leq - R(\bb,a) + \sum_{\bb' \in \mathcal{B}} \pb(\bb'| \bb,a), \nonumber \\
&& \hspace{1.2cm} \mbox{for all } \bb \in \mathcal{B}, \ a \in \mathcal{A}, \ \psi \in \mathbb{R}, \ \phi(\cdot) \in F^{\mu}.  \nonumber
\end{eqnarray}
 We have the following result:
\begin{thm}{\textbf{[Strong duality]}}
\label{thm_strong_duality_uncons}
Assume A.1-2, then 
\begin{enumerate}[(i)]
    \item There exists an average optimal stationary policy;
    \item The LP \eqref{Eqn_dual_lp_uncons} is solvable and there is no duality gap. 
\end{enumerate} 
\end{thm}
\textbf{Proof:} From Theorem \ref{Thm_two_parts}, the  belief MDP has a uniform Lyapunov function. As the immediate rewards are bounded and by the definition of the dual \eqref{Eqn_dual_lp_uncons} itself, $\phi(\cdot)$ is restricted to linear space $F^{\mu}$, then using theorem \cite[ Theorem  12.4]{altman2021constrained}, the dual \eqref{Eqn_dual_lp_uncons} is feasible and has no duality gap. \eop

\subsection{Constrained IOMDPs}
 In the previous subsection \ref{sec_uncons_iomdp}, we discussed the unconstrained belief MDP with average reward criterion. We establish strong duality for the same in Theorem \ref{thm_strong_duality_uncons}. In this subsection, our aim is to prove the strong duality for the constrained belief MDP. Specifically, we consider the belief MDP problem formulated in \eqref{Eqn_average_reward_cons_belief}.    Then, one can define the belief states as in the previous section and the LP that can potentially solve the constrained IOMDP problem \eqref{Eqn_average_reward_cons_belief}:
\begin{eqnarray}
\label{Eqn_primal_lp_cons}    
\inf_{ {\bf x} = \{x(\bb, a)\}_{\bb, a} } && \sum_{\bb,a}   -R(\bb,a) x(\bb,a)\ \  \mbox{subject to }   \\
&&  \sum_{\bb,a} C(\bb, a) x(\bb,a) \le B,   \nonumber \\
&& \sum_a x(\bb',a) = \sum_{\bb,a} Q(\bb' | \bb, a) x(\bb,a) , \ \forall \bb'   \hspace{4mm}\label{Eqn_constraint_second} \\
&& \sum_{\bb,a} x(\bb,a) = 1, \  x(\bb,a) \ge 0, \ \forall \bb, a,  \nonumber 
\end{eqnarray}
where $\bb \in \mathcal{B}$ represents the belief state and  
\begin{equation}
  C(\bb, a):=  \sum_{s \in \s} \bb(s) c(s,a),
  \label{Eqn_belief_cons_cost}
\end{equation}
which is bounded since $c(\cdot, \cdot)$ is bounded. 
 Now following the same procedure as in the previous subsection \ref{sec_uncons_iomdp}, the  dual problem  of \eqref{Eqn_primal_lp_cons} is given by  
\begin{eqnarray}
\label{Eqn_dual_lp_cons}    
 \sup_{\psi, \phi, \lambda \ge 0}  \hspace{-2mm}&\hspace{-2mm}&\hspace{-2mm} \left(\psi - \lambda B \right)\quad  \mbox{subject to}  \\
&\hspace{-2mm}&\hspace{-2mm}\phi(\bb)+ \psi \leq - R(\bb,a) + \lambda C(\bb,a) +\sum_{\bb' \in \mathcal{B}} \pb(\bb'| \bb,a), \nonumber\\
&& \hspace{1.3cm} \mbox{for all } \bb \in \mathcal{B}, a \in \mathcal{A}, \psi \in \mathbb{R}, \phi(\cdot) \in F^{\mu}.   \nonumber
\end{eqnarray}
 We have the following result:
\begin{thm}{\textbf{[Strong duality for constrained IOMDPs]}} 
\label{thm_strong_duality_cons}
Assume A.1-2, then:
\begin{enumerate}[(i)]
    \item There exists an average optimal stationary policy;
    \item The LP \eqref{Eqn_dual_lp_cons} is solvable and there is no duality gap. 
\end{enumerate} 
\end{thm}
\textbf{Proof:}
Here, constraint coefficient $C(\cdot, \cdot)$ is also  bounded.  Thus, the result is followed using same arguments as in proof of Theorem \ref{thm_strong_duality_uncons}, only the  last step is replaced by \cite[ Theorem 12.9]{altman2021constrained}. \eop

\section{Wireless Channel Application}
\label{sec_example}
Consider a scenario where a controller relies on a remote sensor to perceive the state information of a dynamic process. However, the controller and the sensor are physically separated and communicate over an unreliable wireless channel. Due to environmental interference and channel fading, the transmissions over this wireless medium may occasionally fail. Consequently, the state information observed by the sensor is not always successfully delivered to the controller. To address  such challenges, the IOMDP framework of  section \ref{sec_model} is very useful; the framework  provides a structured approach to selecting optimal actions based on  intermittent state information.  
The basic idea is to choose an optimal action based on the last observed state and the age of the information (age is technically defined as  the time gap after the last observation, see e.g., \cite{kavitha2021controlling, kaul2012real}). 

The dynamical system   takes values from a finite state space, denoted as $\mathcal{S} = \{1, 2, \dots, k\}$, where $1$ represents the worst system   state and $k$ represents the best (in terms of deriving the utilities, as explained  immediately below). The intermediate states are progressively better.   The system state evolves on its own over time according to a Markov process with a fixed probability transition model, given by $P_W:= \{p_W(s'|s)\}_{s, s'}$, which is beyond the controller’s influence. To be in accordance with {A.2}, we assume that $P_W$  is  recurrent.  

The controller has to interact with the above system on a regular basis by spending some energy. If it interacts with the energy level $a$ when the system is in state $s$, then it obtains a utility equal to $r(s,a)$. In addition, it has a constraint on the average energy spent and the instantaneous cost of the energy spent is given by $c(a)$. Let $B< \infty$ represent the time-average energy constraint. 


Hypothetically,  say  the controller can observe the  dynamical system  completely.  Then  
it would have been  interested in the following constrained MDP problem:
\begin{eqnarray}
\sup_\pi \hspace{-3mm}&\hspace{-3mm}&\hspace{-3mm} \left(\lim_{T \to \infty} \frac{1}{T} \sum_{t=1}^T \E[ r(s_t, a_t)]  \right), 
\mbox{ subject to }  \nonumber \\
 && \lim_{T \to \infty} \frac{1}{T} \sum_{t=1}^T \E[ c( a_t)]  \le B.  
\end{eqnarray}
However,  the controller  can observe the states only intermittently (recall the information is transferred over wireless channel); hence  in reality, it has to choose  actions based on the intermittent observations or on the  belief states as in section \ref{sec_model}. Thus, we now describe the relevant IOMDP. 
 
As the dynamical system is evolving on its own,  the transition from  the state $s$  to state $s'$ equals  $p_W(s'|s)$, which is independent of the  energy spent. 
Note that, in this application with action-independent transitions, the belief state $\bb$ can be represented by two components $(s,\eta)$, where $s \in \cal S$ is the last observed state and $\eta \in {\mathbb{N}_0}$ is the age of the state information. 
 
The transition probabilities for the belief MDP are:
\begin{equation}
\label{eqn_wireless_tpm}
\pb \left(\bb_t \mid \bb_{t-1}\right) :=   \left\{\begin{array}{ll}
1-\rho, & \text {if } \bb_t=P_W^{\top} \bb_{t-1},  \\
\rho\left[P_W^{\top} \bb_{t-1}\right]_i, & \text {if } \bb_t=e_i, \ i \in \mathcal{S},  \\
0, & \text {otherwise}.
\end{array}\right. 
\end{equation}
The reward and cost functions $R(\cdot, \cdot)$ and  $C(\cdot, \cdot)$ are as in \eqref{Eqn_belief_reward} and \eqref{Eqn_belief_cons_cost}. Note that the age $\eta$ is geometric random variable with parameter $\rho$.  Now define, $\Bar{x}(\bb) := \sum_a x(\bb, a)$.  The   constraint  \eqref{Eqn_constraint_second} for belief chain \eqref{eqn_wireless_tpm} modifies to:
$$
\bar{x}(\bb') = \sum_{\bb} \pb(\bb' | \bb ) \bar{x}(\bb), \ \forall \bb'. 
$$
 The  solution of the above system of equations  with further constraints $\sum_{\bb} \bar{x}(\bb) = 1, \  \bar{x}(\bb) \ge 0,$ for all $\bb$, is the  stationary distribution, call it $\bf{\nu} (\cdot)$ (recall by Theorem \ref{Thm_two_parts}.$(i)$, the belief MDP is unichain).   By Theorem \ref{thm_strong_duality_cons},  the strong duality holds for this wireless IOMDP problem
  and  hence the optimal policy can be derived by solving  the following  LP ({\it c.f.} \eqref{Eqn_primal_lp_cons}): 
\begin{eqnarray}
\label{Eqn_primal_lp_wireless}    
\sup_{ {\bf x} = \{x(\bb, a)\}_{\bb, a} } && \sum_{\bb,a}  x(\bb,a) R(\bb,a)\ \  \mbox{subject to } \nonumber \\
&&  \sum_{\bb,a} C(a) x(\bb,a) \le B   \nonumber \\
&& \sum_{a}  x(\bb,a) = \nu(\bb), \  x(\bb,a) \ge 0, \forall \bb, a. \ \ 
 \end{eqnarray}
Let $\gamma $ represent   the stationary distribution (SD) of the original Markov  chain with transition probability matrix $P_W$. Then, the SD $\nu$ of the belief Markov chain is given by:
$$
\nu(\bb) = {\gamma(s)} \rho (1-\rho)^{\eta}    \mbox{ for any  } \bb = (s, \eta).  
$$
We next move towards solving LP \eqref{Eqn_primal_lp_wireless} numerically.

\section{Numerical Results}
\label{sec_numerical}

In this numerical study, we approximately solve the LP \eqref{Eqn_primal_lp_wireless} and find a near-optimal policy for the application with wireless observations. 
The aim here is to understand the properties of the near-optimal policy related to various system parameters; with the more specific focus being on age dependency.

By established strong duality,  one can use LP \eqref{Eqn_primal_lp_wireless} to find the optimal policy. 
However,  the dual LP is infinite dimensional, and hence one needs to  truncate   \eqref{Eqn_primal_lp_wireless} appropriately to numerically derive  a `good' or near-optimal policy. 
The age variables are geometrically distributed, and hence the natural choice  is to truncate all the belief states with large age components.  Thus, we consider the $K$-truncated approximate LP, 
with  belief states  spanning over $\{ (s, \eta) : s \in {\cal S}, \eta \leq K\} $. 
We leave formal justification of this approximation that can be done with the help of \cite{demirci2024average,hernandez2012adaptive,hernandez2012further} as a topic for future research. 


We consider the following  set the parameters: $\mathcal{S} = \{ 1, 2\}$, $\mathcal{A} = \{ a_1, a_2\}$, $C(a) = (2+a)^2$, $B=10.4$, rewards  and transition probabilities are (respectively)
\begin{eqnarray*}
    [r(s,a)]_{s,a} = \begin{bmatrix}
  0 & 1    \\
1 & 4
\end{bmatrix} \mbox{ and } P_W = \begin{bmatrix}
  0.7 & 0.3    \\
0.1 & 0.9
\end{bmatrix}. 
\end{eqnarray*}
The idea here is to consider a skewed system, where state~$2$ is the first choice with perfect information and limited energy budget, 
and then understand the age-dependent characteristics of the optimal policy. 
Further, from the transition matrix $P_W$, the system has higher tendency to  prolong in state~$2$ than in state~$1$. This aspect can further accentuate the age dependencies and we proceed immediately towards numerically studying the same. 
In the examples of this section, we  consider $\rho \le 0.6$ and hence set $K = 10$; this provides sufficiently good approximation (observe that $0.6^{10}  \approx 0.006$, which is sufficiently small). 

\subsection{Study with respect to age}
  In Tables \ref{table_state1}-\ref{table_state2}, we vary the observation probability $\rho$ and study the corresponding  variations in  the (age-dependent) optimal policy.  In Table \ref{table_state1}, we tabulate the probability  of choosing lower energy $a=a_1$   at different ages and with state $1$ (for brevity, we denote  $\pi^*(s=1, \eta)$ by $d^*_\eta$). As seen from Table \ref{table_state1}, at  $\rho = 0.1$,  it is optimal to choose high  energy for higher ages; observe $d^*_\eta = 0$   for higher ages  $\eta$, when $\rho = 0.1$; this is probably due to the fact that the belief of being in state $2$ increases with age.

\medskip
  
  \begin{table}[h!]
\centering
\begin{tabular}{|l|l|l|l|l|l|l|l|l|l|l|}
\hline
$\bf{\rho}$ &  \bf{$\eta$=0} &\bf{$\eta$=1}  & \bf{$\eta$=2}  &  \bf{$\eta$=3}  &\bf{$\eta$=4}   \\ \hline
0.1   & 1    & 1   & 1    & 0.8208     & 0       \\ \hline
0.2   & 1    & 1   & 1    & 1         & 1       \\ \hline
0.3   & 1    & 1   & 1    & 1         & 1       \\ \hline
0.4   & 1    & 1   & 1    & 1         & 1       \\ \hline
0.5   & 1    & 1   & 1    & 1         & 1       \\ \hline
0.6   & 1    & 1   & 1    & 1         & 1       \\ \hline
\end{tabular}
\caption{State $1$: Optimal policy (probability of low energy {\normalsize{$a_1$}})}
\label{table_state1}
\end{table} 
\vspace{-4mm}
\begin{table}[h!]
\centering
\begin{tabular}{|l|l|l|l|l|l|l|l|l|l|l|}
\hline
$\bf{\rho}$ &  \bf{$\eta$=0} &\bf{$\eta$=1}  & \bf{$\eta$=2}  &  \bf{$\eta$=3}  &\bf{$\eta$=4}   \\  \hline
0.1   & 0       & 0     & 0      & 0         & 0       \\ \hline
0.2   & 0       & 0     & 0.5786  & 1         & 1       \\ \hline
0.3   & 0       & 0.9973 & 1      & 1         & 1       \\ \hline
0.4   & 0.3178   & 1     & 1      & 1         & 1       \\ \hline
0.5   & 0.4650   & 1     & 1      & 1         & 1       \\ \hline
0.6   & 0.5554   & 1     & 1      & 1         & 1       \\ \hline
\end{tabular}
\caption{State $2$: Optimal policy (probability of low energy {\normalsize{$a_1$}})}
\label{table_state2}
\end{table}
   More interestingly, this phenomenon is not observed at higher $\rho = 0.6$, $d^*_\eta = 1$ for all $\eta$; with larger $\rho$, the chances of observing the true state increase; in other words, $P(s=2, \eta = 0)= \nu(2, 0) $ increases and this coupled with the (same) energy constraint could be the cause; basically, with larger $\rho$ high energy $a=a_2$ is chosen for states with more certainty.



In Table \ref{table_state2}, we tabulate the optimal policy corresponding to state~$2$. 
The observations again  mirror similar aspects.  A notable observation is that with $\rho = 0.6$, at state $2$ and age $\eta = 0$ (i.e., with perfect information of being in `optimal' state) the probability of choosing high energy $a=a_2$ is set at the maximum possible value that  satisfies  the energy budget constraint. 

In summary,
with larger $\rho = 0.6$, the high-energy action ($a=a_2$) is chosen only for belief states with more certainty of being in state~$2$---the optimal probability of spending low energy $a_1$ is one, for state $2$ with $\eta \ge 1$ and at all ages with state $1$. However, when $\rho$ decreases, the probability of choosing the low energy option $a_1$ decreases with high ages in state $1$ and increases with high ages in state $2$, where the belief of being in state $2$ is higher.


    \subsection{Study with respect to $\rho$}
In Tables \ref{table_state1}-\ref{table_state2}, we also analyze the optimal policy dependency on the observation probability $\rho$. From the tables, when the age $\eta$ is large,  the probability of spending low energy $d^*_{\eta}$ increases with  observation probability $\rho$. This is again due to energy budget constraint---as $\rho$ increases, it is optimal to utilize high energy action  $a_2$ on belief states with higher certainty of being in state $2$.

\section{Conclusion}
\label{conclusion}

This paper studies the average reward intermittently observable MDPs. We assumed that the state is observed intermittently with a fixed probability. We first formulated the IOMDP as belief MDP, which has countable state space. We proved that the belief Markov chain is unichain, when the original MDP is recurrent.  Under the same recurrence assumption, we have established that the  strong duality holds for both unconstrained and constrained belief MDPs. 

Finally, a wireless example with observation erasures is provided,  where the optimal policy has to choose the energy to be spent on utilizing the system depending upon the age of the information and the last observed state. 
We solve the problem approximately and numerically using our LP approach and make some interesting observations about the structure of the near-optimal policy.

\section{Future Work}
\label{sec_future_work}
As a first direction, we would like to extend the present results by considering a multi-chain model in place of Assumption~A.2. Then, it is interesting to investigate the sensitivity of the proposed framework with respect to variations in system parameters and uncertainty factors that may arise in practical applications. Another line of extension will focus on expanding the numerical experiments to multi-state settings, thereby enabling a more comprehensive evaluation of the policy structure. In addition, we aim to study and compare the performance of different truncation strategies, including adaptive approaches where the truncation level 
$K$ is dynamically adjusted.

\vspace{-2mm}
\section{Appendix }\label{sec_app_a}
\textbf{Proof of Theorem \ref{Thm_two_parts}:} \textbf{Proof of part (i):} Fix any SMR policy $\pi$. Then the following is true. 

\textit{Existence of single recurrent class:} To begin with, observe that, under  A.2, the pure states, i.e.,  $\bb \in \mathcal{E} = \{e_i : i \in \mathcal{S}\}$ are all  recurrent. Thus, there exists at least one closed recurrent class.  We now prove by contradiction that there exists unique recurrent class. Suppose there exist two closed recurrent classes, say  ${\cal{C}}_1$ and ${\cal{C}}_2$, both containing pure belief states.
Without loss of generality, assume   $\{ e_i: i \leq n < |\mathcal{S}|\} \subset {\cal{C}}_1$ and $\{ e_i: i > n \} \subset {\cal{C}}_2$. 

By assumption~{A.2}, there exists an $n>0$ such that  $P_a^n (i,j) > 0$ for all $i,j \in \mathcal{S}$ and all $a \in \mathcal{A}$. Thus, the belief transition probability from \eqref{eqn_wireless_tpm}, satisfies 
$$
\pb^n(e_j \mid e_i, a) > \rho \min_{a'} P_{a'}^n (i,j) > 0, \mbox{ for all } i, j \in \mathcal{S} \mbox{ and }  a \in \mathcal{A}.
$$ 
Consequently, there exists a strictly positive probability of transitioning from ${\cal{C}}_1$ to ${\cal{C}}_2$ and vice versa. This contradicts the definition of closed recurrent classes.  It is clear  (for example) that  ${\cal{C}}_2$ is not recurrent if ${\cal{E}} \subset {\cal{C}}_1$. 
 
Let $\mathcal{C}$ be the unique recurrent class. Next, we establish that any state outside $\mathcal{C}$ must be transient. Since $\mathcal{E} \subset \mathcal{C}$, any belief state $\bb \notin \mathcal{C}$ must satisfy $\bb \notin \mathcal{E}$. From the belief transition probability  \eqref{eqn_belief_trans_probability}, there is a strictly positive probability of transitioning from such a belief state $\bb$ to some $e_i \in \mathcal{C}$. This implies that $\bb$ is a transient state. Since the time between two consecutive state observations is geometrically distributed, the absorption expected time to the single class is finite. Hence, the assumption $B1$ from \cite[Chapter 11]{altman2021constrained} is satisfied. 
 

{\textit{Positive recurrence:}} Now, we consider the closed recurrent class $\mathcal{C}$ and focus on the belief chain evolving only over $\mathcal{C}$ as state space. Take set $\mathcal{E}$ and consider the function $\mu: \mathcal{C} \to [1,\infty)$ such that
\begin{equation}
\label{mufosters}
\mu(\bb) =   \left\{\begin{array}{ll}
{2}, & \bb \in \mathcal{C}\setminus \mathcal{E},\\
1, & \bb \in \mathcal{E}. 
\end{array}\right.    
\end{equation}
As $\mathcal{S}$ is finite set, $\mathcal{E}$ is also finite. Using Foster's criteria \cite[Theorem 2.6.4]{menshikov2016non}, for $\bb \in \mathcal{C}\setminus \mathcal{E}$
\vspace{-1mm}
{\small
\begin{eqnarray*}
    \mathbb{E}[ \mu(\bb_{t+1}) - \mu(\bb_t) | \bb_t=\bb] \hspace{-2.5mm}&\hspace{-2.5mm}=\hspace{-2.5mm}&\hspace{-2.5mm} \sum_{\bb' \in \mathcal{C}} \pb(\bb'|\bb, \pi(\bb)) (\mu(\bb')-\mu(\bb))  \\
    &=& - \sum_{\bb' \in \mathcal{E}}  \pb(\bb'|\bb,  \pi(\bb)).  \\
\end{eqnarray*}}
Also,
$\mathbb{E}[ \mu(\bb_{t+1})| \bb_t=\bb] = \sum_{\bb' \in  \mathcal{C}} \pb(\bb'|\bb,  \pi(\bb)) \mu(\bb') \leq 2,$  for all $\bb \in  \mathcal{E}$. 
As belief chain is irreducible over $\mathcal{C}$, using \cite[Theorem 2.6.4]{menshikov2016non}, it is positive recurrent. 

 \textbf{Proof of part (ii):} 
Consider the function $\mu: \mathcal{B} \to [1, \infty)$, such that 
\begin{equation*}
\mu(\bb) =   \left\{\begin{array}{ll}
{1}, & \bb \in \mathcal{B} \setminus \mathcal{E}, \nonumber\\
2, & \bb \in \mathcal{E}. 
\end{array}\right.   
\end{equation*}
 
Now, for contraction mapping, we consider for any $\bb \in \mathcal{B}$, 
\begin{eqnarray*}
   \sum_{\bb' \notin \mathcal{E}} \pb(\bb'|\bb,  a) \mu(\bb') =    (1- \rho) \leq (1- \rho) \mu(b). 
\end{eqnarray*}
Here, $(1 - \rho) \in [0,1)$. Also, 1) $\langle  \beta, \mu\rangle <\infty$ as both $\beta$ (initial distribution) and $\mu$ are bounded. 2) As $\mathcal{A}$ set is finite, then under discrete topology, the following condition is also satisfied, as $a_n \to a$
$$
\lim_{n \to \infty} \sum_{\bb' \in \mathcal{B}} |\pb (\bb'| \bb, a_n)- \pb (\bb'| \bb, a)| \mu(\bb') =0.
$$
3) Now, from the boundness of $R(\cdot, \cdot)$ and the boundness of $\mu(\cdot)$ (see \eqref{mufosters}), we conclude that $R(\cdot, \cdot)$ is $\mu$-bounded. Thus, all the conditions for the contracting MDP are satisfied.   

\textit{Proof of uniform Lyapunov function}: Now, by   \cite[ Lemma 11.5]{altman2021constrained}, the following conditions are satisfied (equivalent to $B2$ and $B3$ in \cite[Chapter 11]{altman2021constrained}):
\begin{enumerate}
    \item Given a policy $\pi$ and initial distribution $\beta$, the expected occupation measures $\{f^t(\beta, \pi) \}_t$ are tight, where 
    $$
f^t(\beta, \pi)(\bb, a) := \frac{1}{t} \sum_{s=1}^t \pb_{\beta}^{\pi}(\bb_s = \bb, a_s = a),  
    $$
    with $\pb_{\beta}^{\pi}$ is probability generated by initial distribution $\beta$ and policy $\pi$. 
    \item Given a policy $\pi$ and initial distribution $\beta$, the expected occupation measures $\{f^t(\beta, \pi) \}_t$ are integrable with respect to the absolute value of the immediate reward $R$, uniformly in $t$. 
\end{enumerate}
    As  both  the above conditions are satisfied, by  \cite[Theorem  11.14]{altman2021constrained}, the belief MDP has a uniform Lyapunov function. \eop

\bibliographystyle{IEEEtran}
\bibliography{ref}

\begin{thebibliography}{10}
\providecommand{\url}[1]{#1}
\csname url@samestyle\endcsname
\providecommand{\newblock}{\relax}
\providecommand{\bibinfo}[2]{#2}
\providecommand{\BIBentrySTDinterwordspacing}{\spaceskip=0pt\relax}
\providecommand{\BIBentryALTinterwordstretchfactor}{4}
\providecommand{\BIBentryALTinterwordspacing}{\spaceskip=\fontdimen2\font plus
\BIBentryALTinterwordstretchfactor\fontdimen3\font minus \fontdimen4\font\relax}
\providecommand{\BIBforeignlanguage}[2]{{%
\expandafter\ifx\csname l@#1\endcsname\relax
\typeout{** WARNING: IEEEtran.bst: No hyphenation pattern has been}%
\typeout{** loaded for the language `#1'. Using the pattern for}%
\typeout{** the default language instead.}%
\else
\language=\csname l@#1\endcsname
\fi
#2}}
\providecommand{\BIBdecl}{\relax}
\BIBdecl

\bibitem{puterman2014markov}
M.~L. Puterman, \emph{Markov decision processes: {D}iscrete stochastic dynamic programming}.\hskip 1em plus 0.5em minus 0.4em\relax John Wiley \& Sons, 2014.

\bibitem{krishnamurthy2016partially}
V.~Krishnamurthy, \emph{Partially observed {M}arkov decision processes}.\hskip 1em plus 0.5em minus 0.4em\relax Cambridge University Press, 2016.

\bibitem{borkar2003dynamic}
V.~S. Borkar, ``Average cost dynamic programming equations for controlled {M}arkov chains with partial observations,'' \emph{SIAM Journal on Control and Optimization}, vol.~39, no.~3, pp. 673--681, 2000.

\bibitem{borkar2007dynamic}
------, ``Dynamic programming for ergodic control of {M}arkov chains under partial observations: {A} correction,'' \emph{SIAM Journal on Control and Optimization}, vol.~45, no.~6, pp. 2299--2304, 2007.

\bibitem{demirci2024average}
Y.~E. Demirci, A.~D. Kara, and S.~Y{\"u}ksel, ``Average cost optimality of partially observed {MDPs}: {C}ontraction of nonlinear filters and existence of optimal solutions and approximations,'' \emph{SIAM Journal on Control and Optimization}, vol.~62, no.~6, pp. 2859--2883, 2024.

\bibitem{hsu2006existence}
S.-P. Hsu, D.-M. Chuang, and A.~Arapostathis, ``On the existence of stationary optimal policies for partially observed {MDPs} under the long-run average cost criterion,'' \emph{Systems \& Control Letters}, 2006.

\bibitem{kara2022near}
A.~Kara and S.~Y{\"u}ksel, ``Near optimality of finite memory feedback policies in partially observed {M}arkov decision processes,'' \emph{Journal of Machine Learning Research}, vol.~23, no.~11, pp. 1--46, 2022.

\bibitem{platzman1980optimal}
L.~K. Platzman, ``Optimal infinite-horizon undiscounted control of finite probabilistic systems,'' \emph{SIAM Journal on Control and Optimization}, vol.~18, no.~4, pp. 362--380, 1980.

\bibitem{runggaldier1993nearly}
W.~J. Runggaldier and L.~Stettner, ``Nearly optimal controls for stochastic ergodic problems with partial observation,'' \emph{SIAM Journal on Control and Optimization}, vol.~31, no.~1, pp. 180--218, 1993.

\bibitem{chen2023intermittently}
G.~Chen and S.-C. Liew, ``Intermittently observable {M}arkov decision processes,'' \emph{arXiv preprint arXiv:2302.11761}, 2023.

\bibitem{khan2023cooperative}
N.~Khan and V.~Subramanian, ``Cooperative multi-agent constrained {POMDPs}: {S}trong duality and primal-dual reinforcement learning with approximate information states,'' \emph{arXiv preprint arXiv:2307.16536}.

\bibitem{shuman2010measurement}
D.~Shuman and et~al, ``Measurement scheduling for soil moisture sensing: From physical models to optimal control,'' \emph{Proceedings of the IEEE}, vol.~98, no.~11, pp. 1918--1933, 2010.

\bibitem{bertsekas1996stochastic}
D.~Bertsekas and S.~E. Shreve, \emph{Stochastic optimal control: {T}he discrete-time case}.\hskip 1em plus 0.5em minus 0.4em\relax Athena Scientific, 1996, vol.~5.

\bibitem{altman2021constrained}
E.~Altman, \emph{Constrained {M}arkov decision processes}.\hskip 1em plus 0.5em minus 0.4em\relax CRC, 1999.

\bibitem{lasserre1994average}
J.~B. Lasserre, ``Average optimal stationary policies and linear programming in countable space {M}arkov decision processes,'' \emph{Journal of Mathematical Analysis and Applications}, vol. 183, no.~1, pp. 233--249.

\bibitem{lim2023optimality}
M.~H. Lim and et~al, ``Optimality guarantees for particle belief approximation of {POMDPs},'' \emph{Journal of Artificial Intelligence Research}, vol.~77, pp. 1591--1636, 2023.

\bibitem{ahmadi2020control}
M.~Ahmadi, N.~Jansen, B.~Wu, and U.~Topcu, ``Control theory meets {POMDPs}: {A} hybrid systems approach,'' \emph{IEEE Transactions on Automatic Control}, vol.~66, no.~11, pp. 5191--5204, 2020.

\bibitem{akbarzadeh2024two}
N.~Akbarzadeh and A.~Mahajan, ``Two families of indexable partially observable restless bandits and {W}hittle index computation,'' \emph{Performance Evaluation}, vol. 163, p. 102394, 2024.

\bibitem{kavitha2021controlling}
V.~Kavitha and E.~Altman, ``Controlling packet drops to improve freshness of information,'' in \emph{Proceedings of NetGCooP 2020}.

\bibitem{kaul2012real}
S.~Kaul, R.~Yates, and M.~Gruteser, ``Real-time status: {H}ow often should one update?'' in \emph{Proceedings IEEE INFOCOM 2012}.

\bibitem{hernandez2012adaptive}
O.~Hern{\'a}ndez-Lerma, \emph{Adaptive {M}arkov control processes}.\hskip 1em plus 0.5em minus 0.4em\relax Springer, 2012.

\bibitem{hernandez2012further}
O.~Hern{\'a}ndez-Lerma and J.~B. Lasserre, \emph{Further topics on discrete-time {M}arkov control processes}.\hskip 1em plus 0.5em minus 0.4em\relax Springer, 2012.

\bibitem{menshikov2016non}
M.~Menshikov, S.~Popov, and A.~Wade, \emph{Non-homogeneous random walks: {L}yapunov function methods for near-critical stochastic systems}.\hskip 1em plus 0.5em minus 0.4em\relax Cambridge University Press, 2016.

\end{thebibliography}


\begin{thebibliography}{9}
\bibitem{puterman1990markov}Puterman, M. L. (1990). Markov decision processes. Handbooks in operations research and management science, 2, 331-434.

\bibitem{krishnamurthy2016partially}Krishnamurthy, V. (2016). Partially observed Markov decision processes. Cambridge university press.

\bibitem{borkar2000average}Borkar, V. S. (2000). Average cost dynamic programming equations for controlled Markov chains with partial observations. SIAM Journal on Control and Optimization, 39(3), 673-681.

\bibitem{borkar2007dynamic}Borkar, V. S. (2007). Dynamic programming for ergodic control of Markov chains under partial observations: A correction. SIAM Journal on Control and Optimization, 45(6), 2299-2304.

\bibitem{demirci2024average}Demirci, Y. E., Kara, A. D., and  Yüksel, S. (2024). Average Cost Optimality of Partially Observed MDPs: Contraction of Nonlinear Filters and Existence of Optimal Solutions and Approximations. SIAM Journal on Control and Optimization, 62(6), 2859-2883.

\bibitem{hsu2006existence}Hsu, S. P., Chuang, D. M.,  and  Arapostathis, A. (2006). On the existence of stationary optimal policies for partially observed MDPs under the long-run average cost criterion. Systems  and control letters, 55(2), 165-173.

\bibitem{kara2022near} Kara, A.,  and Yuksel, S. (2022). Near optimality of finite memory feedback policies in partially observed Markov decision processes. Journal of Machine Learning Research, 23(11), 1-46.

\bibitem{platzman1980optimal} Platzman, L. K. (1980). Optimal infinite-horizon undiscounted control of finite probabilistic systems. SIAM Journal on Control and Optimization, 18(4), 362-380.

\bibitem{runggaldier1993nearly} Runggaldier, W. J., and  Stettner, L. (1993). Nearly optimal controls for stochastic ergodic problems with partial observation. SIAM journal on control and optimization, 31(1), 180-218.

\bibitem{chen2023intermittently} Chen, G., and Liew, S. C. (2023). Intermittently Observable Markov Decision Processes. arXiv preprint arXiv:2302.11761.

\bibitem{khan2023cooperative} Khan, N., and Subramanian, V. (2023). Cooperative multi-agent constrained POMDPs: Strong duality and primal-dual reinforcement learning with approximate information states. arXiv preprint arXiv:2307.16536.

\bibitem{bertsekas1996stochastic} Bertsekas, D., and Shreve, S. E. (1996). Stochastic optimal control: the discrete-time case (Vol. 5). Athena Scientific.

\bibitem{altman2021constrained} Altman, E. (2021). Constrained Markov decision processes. Routledge.

\bibitem{lasserre1994average} Lasserre, J. B. (1994). Average optimal stationary policies and linear programming in countable space Markov decision processes. Journal of mathematical analysis and applications, 183(1), 233-249.

\bibitem{lim2023optimality} Lim, M. H., Becker, T. J., Kochenderfer, M. J., Tomlin, C. J., and  Sunberg, Z. N. (2023). Optimality guarantees for particle belief approximation of pomdps. Journal of Artificial Intelligence Research, 77, 1591-1636.

\bibitem{}

\bibitem{}

\bibitem{ahmadi2020control} Ahmadi, M., Jansen, N., Wu, B., and  Topcu, U. (2020). Control theory meets pomdps: A hybrid systems approach. IEEE Transactions on Automatic Control, 66(11), 5191-5204.

\bibitem{kavitha2021controlling} Kavitha, V.,  and  Altman, E. (2021). Controlling packet drops to improve freshness of information. In Network Games, Control and Optimization: 10th International Conference, NetGCooP 2020, France, September 22–24, 2021, Proceedings 10 (pp. 60-77). Springer International Publishing.

\bibitem{kaul2012real} Kaul, S., Yates, R., and Gruteser, M. (2012, March). Real-time status: How often should one update?. In 2012 Proceedings IEEE INFOCOM (pp. 2731-2735). IEEE.

\bibitem{hernandez2012adaptive} Hernández-Lerma, O. (2012). Adaptive Markov control processes (Vol. 79). Springer Science and Business Media.

\bibitem{hernandez2012further} Hernández-Lerma, O., and Lasserre, J. B. (2012). Further topics on discrete-time Markov control processes (Vol. 42). Springer Science and Business Media.

\bibitem{menshikov2016non} Menshikov, M., Popov, S., and Wade, A. (2016). Non-homogeneous random walks: Lyapunov function methods for near-critical stochastic systems (Vol. 209). Cambridge University Press.

\end{thebibliography}

{\ignore{
}}

\end{document}